\documentclass[12pt, nonatbib]{elsarticle}

\usepackage{amsmath}
\usepackage{mathscinet}
\usepackage{amssymb}
\usepackage{amsthm}
\usepackage{graphicx}
\usepackage[colorlinks]{hyperref}
\AtBeginDocument{%
  \hypersetup{%
    linkcolor=blue,%
    citecolor=green,%
  }%
}
\usepackage{cleveref}
\usepackage{geometry}
\makeatletter
\let\c@author\relax
\makeatother
\usepackage[backend=bibtex,style=numeric,sorting=nyt,
giveninits=true,isbn=false,url=false]{biblatex}
\renewbibmacro{in:}{\ifentrytype{article}{}{\printtext{\bibstring{in}\intitlepunct}}}
\DeclareFieldFormat*{title}{#1}
\crefname{equation}{}{}
\crefname{figure}{}{}

\newtheorem{theorem}{Theorem}[section]
\newtheorem{lemma}[theorem]{Lemma}

\newtheorem{corollary}[theorem]{Corollary}

\theoremstyle{definition}
\newtheorem{definition}[theorem]{Definition}

\theoremstyle{remark}
\newtheorem{remark}[theorem]{Remark}

\numberwithin{equation}{section}

\journal{~~}

\begin{document}

\begin{frontmatter}

\title{Asymptotic behavior for fully nonlinear elliptic equations in exterior domains\tnoteref{t1}}


\author[rvt1]{Yuanyuan Lian}
\ead{lianyuanyuan.hthk@gmail.com}
\author[rvt2]{Kai Zhang\corref{cor1}}
\ead{zhangkaizfz@gmail.com}

\address[rvt1]{Departamento de An\'{a}lisis Matem\'{a}tico,
Instituto de Matem\'{a}ticas IMAG, Universidad de Granada}
\address[rvt2]{Departamento de Geometr\'{i}a y Topolog\'{i}a,
Instituto de Matem\'{a}ticas IMAG, Universidad de Granada}
\tnotetext[t1]{This research has been financially supported by the Project PID2020-118137GB-I00 funded by MCIN/AEI /10.13039/501100011033.}
\cortext[cor1]{Corresponding author. MR Author ID: \href{https://mathscinet.ams.org/mathscinet/2006/mathscinet/search/author.html?mrauthid=1098004}
{1098004}}
\begin{abstract}
In this paper, we obtain the asymptotic behavior at infinity for viscosity solutions of fully nonlinear elliptic equations in exterior domains. We show that if the solution $u$ grows linearly, there exists a linear polynomial $P$ such that $u-P$ is controlled by fundamental solutions of the Pucci's operators. In addition, with proper ellipticity constants, $u(x)-P(x)\to 0$ as $x\to \infty$ (see \Cref{th1.0}). If $u$ grows quadratically, we obtain similar asymptotic behavior (see \Cref{th1.1}). In this paper, we don't require any smoothness of the fully nonlinear operator.
\end{abstract}

\begin{keyword}
Asymptotic behavior \sep exterior domain\sep Liouville theorem \sep fully nonlinear elliptic equation \sep viscosity solution
\MSC[2020] 35B40, 35B53, 35J15, 35J60, 35D40
\end{keyword}

\end{frontmatter}

\section{Introduction}\label{S1}
In this paper, we investigate the asymptotic behavior at infinity of viscosity solutions of fully nonlinear elliptic equations in an exterior domain:
\begin{equation}\label{e.F1}
F(D^2u)=0~~\mbox{ in}~~\bar B_1^c:=\mathbb{R}^n\backslash \bar B_1,
\end{equation}
where $B_1\subset \mathbb{R}^n$ ($n\geq 2$) is the unit open ball and $F$ is always assumed uniformly elliptic with $0<\lambda\leq \Lambda<+\infty$ (see \Cref{de1.1}).

The classical Liouville theorem states that if $u$ is harmonic in $\mathbb{R}^n$ and bounded below or above, it must be a constant. Similarly, if $u$ growth linearly (or quadratically) at most, it must be a linear (or quadratic) polynomial. For $u$ being only harmonic in an exterior domain (i.e., the complement of a compact set), it is not a polynomial in general. Instead, $u$ tends to a polynomial at infinity. This is the counterpart of the Liouville theorem in exterior domains.

The Liouville theorems in exterior domains have been obtained for many equations, such as fully nonlinear uniformly elliptic equations \cite{MR2663711,MR4038557}, the $p$-Laplace equations \cite{MR2424533}, the minimal surface equation \cite{MR898048}, the maximal surface equation \cite{MR4201294}, the Monge-Amper\`{e} equation \cite{MR3299174,MR1953651,MR1679973}, the special Lagrangian equations \cite{MR4038557} etc.

In \cite{MR2663711}, the asymptotic behavior at infinity has been investigated thoroughly if the solution is bounded above or below (see \Cref{pr2.3}). In this paper, we study the asymptotic behavior under the assumption that the solution grows linearly or quadratically. As far as we know, there is no such research without smoothness assumption on the fully nonlinear operator $F$. There are indeed some important nonsmooth operators, such as the Pucci's operators, the Bellman operators and the Isaacs operators (see \cite[Chapter 2.3]{MR1351007}).

We describe the idea of this paper roughly. Once a solution $u$ is bounded below or above, the asymptotic behavior of $u$ can be completely classified (\cite[Theorem 1.10]{MR2663711} or \Cref{pr2.3}) or controlled (\Cref{le2.3}). This is obtained by the Harnack inequality and the comparison principle basically. In addition, if $u$ grows linearly or quadratically at infinity, there exists a linear or quadratic polynomial $P$ such that $u-P$ is bounded below or above. This is derived by comparing $u$ with a sequence of solutions $w_i$ defined in $B_i$. By letting $i\to\infty$ and the Liouville theorem in the whole space, $w_i$ converges to a polynomial $P$. Hence, $u-P$ is bounded below or above. This idea is highly motivated by \cite{MR4201294} (precisely P. 598--599). By combining above arguments together, we can characterize the asymptotic behavior of $u-P$.

Firs, we introduce some definitions.
\begin{definition}\label{de1.1}
The $F$ is called uniformly elliptic with positive constants $0<\lambda\leq\Lambda<+\infty$ if for any $M,N\in \mathcal{S}^n$,
\begin{equation}\label{e1.0}
\mathcal{M}^-(M,\lambda,\Lambda)\leq F(M+N)-F(N)
\leq \mathcal{M}^+(M,\lambda,\Lambda).
\end{equation}
Here, $\mathcal{S}^n$ denotes the set of $n\times n$ symmetric matrices and $\mathcal{M}^-$, $\mathcal{M}^+$ are Pucci's extremal operators, i.e.,
\begin{equation*}
\mathcal{M}^-(M,\lambda,\Lambda):=\lambda\sum_{\lambda_i>0} \lambda_i+
\Lambda\sum_{\lambda_i<0} \lambda_i,~
\mathcal{M}^+(M,\lambda,\Lambda):=\Lambda\sum_{\lambda_i>0} \lambda_i+
\lambda\sum_{\lambda_i<0} \lambda_i,~
\end{equation*}
where $\lambda_i$ are eigenvalues of $M\in \mathcal{S}^n$.

We write $u\in S(\lambda,\Lambda,0)$ if $\mathcal{M}^-(D^2u,\lambda,\Lambda)\leq 0$ and $\mathcal{M}^+(D^2u,\lambda,\Lambda)\geq 0$ in the viscosity sense.
\end{definition}

We consider viscosity solutions throughout this paper and refer to \cite{MR1351007,MR1118699} for more details of the theory of viscosity solutions. A constant is called universal if it depends only on $n,\lambda$ and $\Lambda$.

The fundamental solutions play a key role in studying the asymptotic behavior of solutions near infinity or a singular point. Roughly speaking, what kinds of fundamental solutions an equations has determine what kinds of asymptotic behavior. Armstrong, Sirakov and Smart \cite[Theorem 1.3]{MR2663711} proved the following existence and uniqueness of fundamental solutions for general fully nonlinear uniformly elliptic equations. We use the same notations and notions as in \cite{MR2663711}. Recall that $F$ is called positively homogeneous of degree $1$ if
\begin{equation*}
F(tM)=tF(M),~\forall ~t>0,~\forall ~M\in \mathcal{S}^n.
\end{equation*}
\begin{lemma}\label{pr1.1}
Suppose that $F$ is positively homogeneous of degree $1$. Then there exists a viscosity solution $\Phi$ of
\begin{equation}\label{e.F2}
F(D^2u)=0~~\mbox{ in}~~\mathbb{R}^n\backslash \{0\}
\end{equation}
such that:\\
(i) $\Phi$ is bounded below in $B_1$ and bounded above in $B_1^c$;\\
(ii) For any $t>0$ and $x\in \mathbb{R}^n\backslash \{0\}$,
\begin{equation}\label{e1.2-1}
\Phi(tx)=t^{-\alpha^{*}}\Phi(x) ~~\mbox{  or  } ~~\Phi(tx)=\Phi(x)-\ln t,
\end{equation}
where $\alpha^{*}\in (-1,+\infty)\backslash \{0\}$ depends only on $n$ and $F$, and $\alpha^*\Phi>0$.

Moreover, any viscosity solution satisfying (i) can be written as: $a\Phi(x)+b$ for some $a>0$ and $b\in \mathbb{R}$.
\end{lemma}

\begin{remark}\label{re1.0}
Since $F$ is uniformly elliptic, $\Phi\in C^{1,\alpha}_{loc}(\mathbb{R}^n\backslash \{0\})$ for some universal constant $0<\alpha<1$ (see \cite[Corollary 5.7]{MR1351007}).
\end{remark}

\begin{remark}\label{re1.3}
The $\alpha^*$ is called the scaling exponent of $F$, which plays a key role in the study of asymptotic behavior of solutions. We set $\alpha^*=0$ for the second case of \cref{e1.2-1}.
\end{remark}

\begin{remark}\label{re1.1}
The (normalized) $\Phi$ with
\begin{equation*}
\min_{x\in \partial B_1} (\text{sign} \alpha^* \cdot \Phi(x))=1~~ (\int_{\partial B_1} \Phi=0~~\mbox{if } \alpha^*=0)
\end{equation*}
is called the upward-pointing fundamental solution of \cref{e.F2}. In this note, $\Phi$ always denotes the upward-pointing fundamental solution.
\end{remark}

~\\

The $\tilde{F}(M):=-F(-M)$ is called the dual operator of $F$ and clearly $\tilde{F}$ is also uniformly elliptic and positively homogeneous of degree $1$. Hence, by \Cref{pr1.1}, there exists the upward-pointing fundamental solution $\tilde{\Phi}$ of
\begin{equation}\label{e1.1-2}
\tilde{F}(D^2u)=0~~\mbox{ in }~~\mathbb{R}^n\backslash \{0\}
\end{equation}
and $\tilde{\Phi}$ satisfies (i) and (ii) in \Cref{pr1.1} for some scaling exponent $\tilde{\alpha}^*\in (-1,+\infty)\backslash \{0\}$. Note that $\alpha^*\neq \tilde{\alpha}^*$ in general.

If $u$ is a solution of \cref{e.F2}, $-u$ will be a solution of \cref{e1.1-2} and vice versa. Thus, $-\tilde{\Phi}$ is also a solution of \cref{e.F2}, which is called the downward-pointing fundamental solution of \cref{e.F2}. Therefore, any operator $F$ possesses exactly a pair of normalized fundamental solutions $(\Phi, -\tilde{\Phi})$: the upward-pointing and the downward-pointing fundamental solutions. Based on the fundamental solutions, we can study the asymptotic behavior of $u$ at infinity.

The Pucci's operators are two special and important fully nonlinear elliptic operators. Note that the dual operator of $\mathcal{M}^+$ is $\mathcal{M}^-$, and vice versa. By direct calculation (see \cite{MR1867613}), the upward-pointing fundamental solutions of $\mathcal{M}^+$ and $\mathcal{M}^-$ are:
\begin{equation*}
E^+(x)=\left\{
  \begin{aligned}
&|x|^{1-(n-1)\lambda/\Lambda}  &&~~\mbox{ if }~~\Lambda/\lambda<(n-1),\\
&-\ln|x|   &&~~\mbox{ if }~~\Lambda/\lambda=(n-1),\\
&-|x|^{1-(n-1)\lambda/\Lambda}  &&~~\mbox{ if }~~\Lambda/\lambda>(n-1)\\
  \end{aligned}
  \right.
\end{equation*}
and
\begin{equation*}
E^-(x)=\left\{
  \begin{aligned}
&-\ln|x|  &&~~\mbox{ if }~~\Lambda/\lambda=1,n=2,\\
&|x|^{1-(n-1)\Lambda/\lambda}  &&~~\mbox{ other cases}\\
  \end{aligned}
  \right.
\end{equation*}
respectively. Correspondingly, the downward-pointing fundamental solutions are
\begin{equation*}
e^+(x)=-E^-(x)=\left\{
  \begin{aligned}
&\ln|x|  &&~~\mbox{ if }~~\Lambda/\lambda=1,n=2,\\
&-|x|^{1-(n-1)\Lambda/\lambda}  &&~~\mbox{ other cases}
  \end{aligned}
  \right.
\end{equation*}
and
\begin{equation*}
e^-(x)=-E^+(x)=\left\{
  \begin{aligned}
&-|x|^{1-(n-1)\lambda/\Lambda}  &&~~\mbox{ if }~~\Lambda/\lambda<(n-1),\\
&\ln|x|    &&~~\mbox{ if }~~\Lambda/\lambda=(n-1),\\
&|x|^{1-(n-1)\lambda/\Lambda}  &&~~\mbox{ if }~~\Lambda/\lambda>(n-1)\\
  \end{aligned}
  \right.
\end{equation*}
respectively.

We use $\alpha^+$ and $\alpha^-$ to denote the scaling exponents of $\mathcal{M}^+$ and $\mathcal{M}^-$ respectively. Then
\begin{equation*}
\frac{(n-1)\lambda}{\Lambda}-1=\alpha^+\leq
\alpha^-=\frac{(n-1)\Lambda}{\lambda}-1.
\end{equation*}
Note that we always have $\alpha^-\geq 0$ since $\Lambda/\lambda\geq 1$ and $n\geq 2$.
Furthermore, for any operator $F$ (see \cite[P.741]{MR2663711}),
\begin{equation*}
\alpha^+\leq \alpha^*(F)\leq \alpha^-,
\end{equation*}
where $\alpha^*(F)$ is the scaling exponent of $F$.

\begin{figure}[h]\label{f.1}
  \centering
  \includegraphics[width=60 mm]{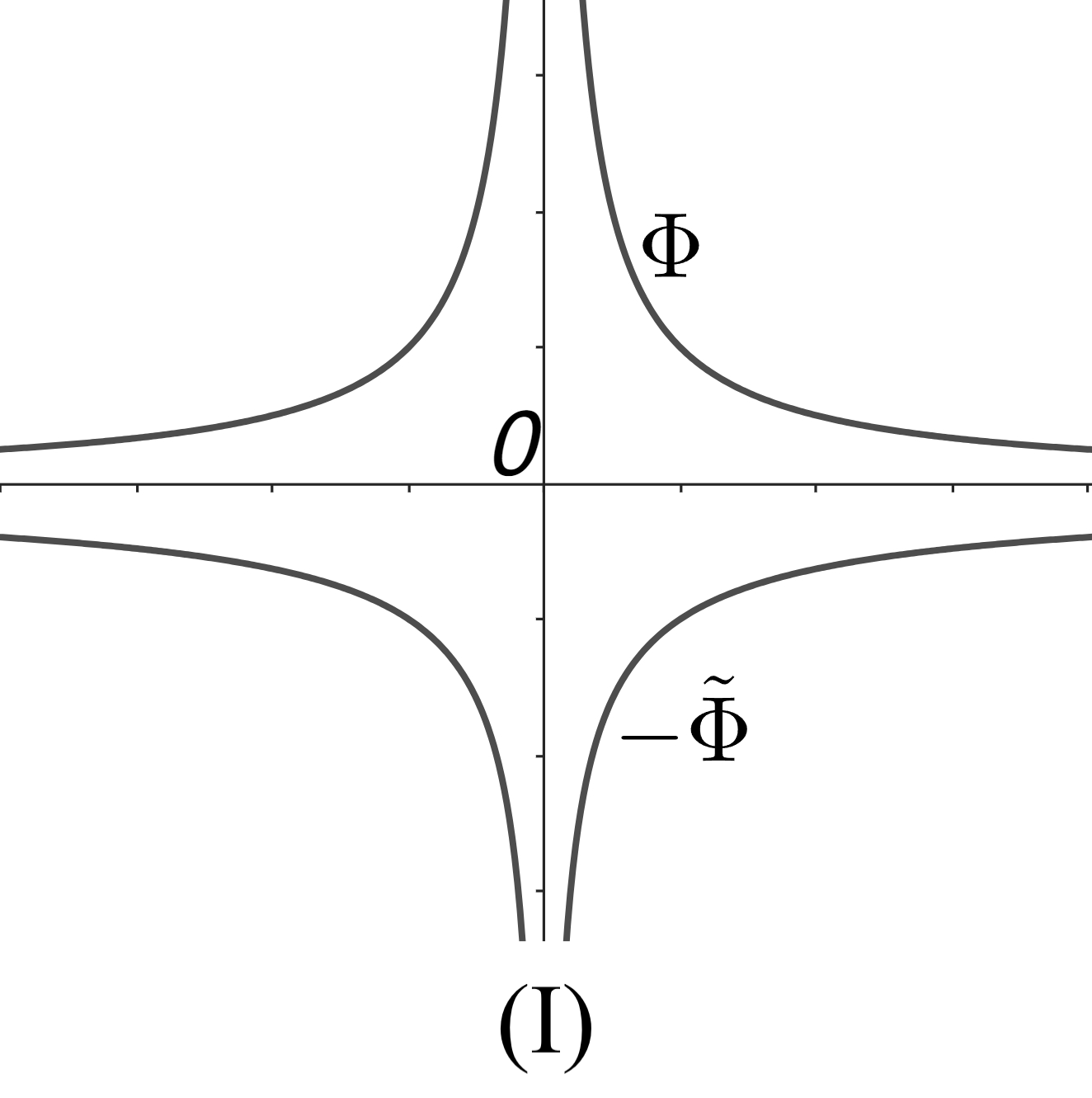}
  \includegraphics[width=60 mm]{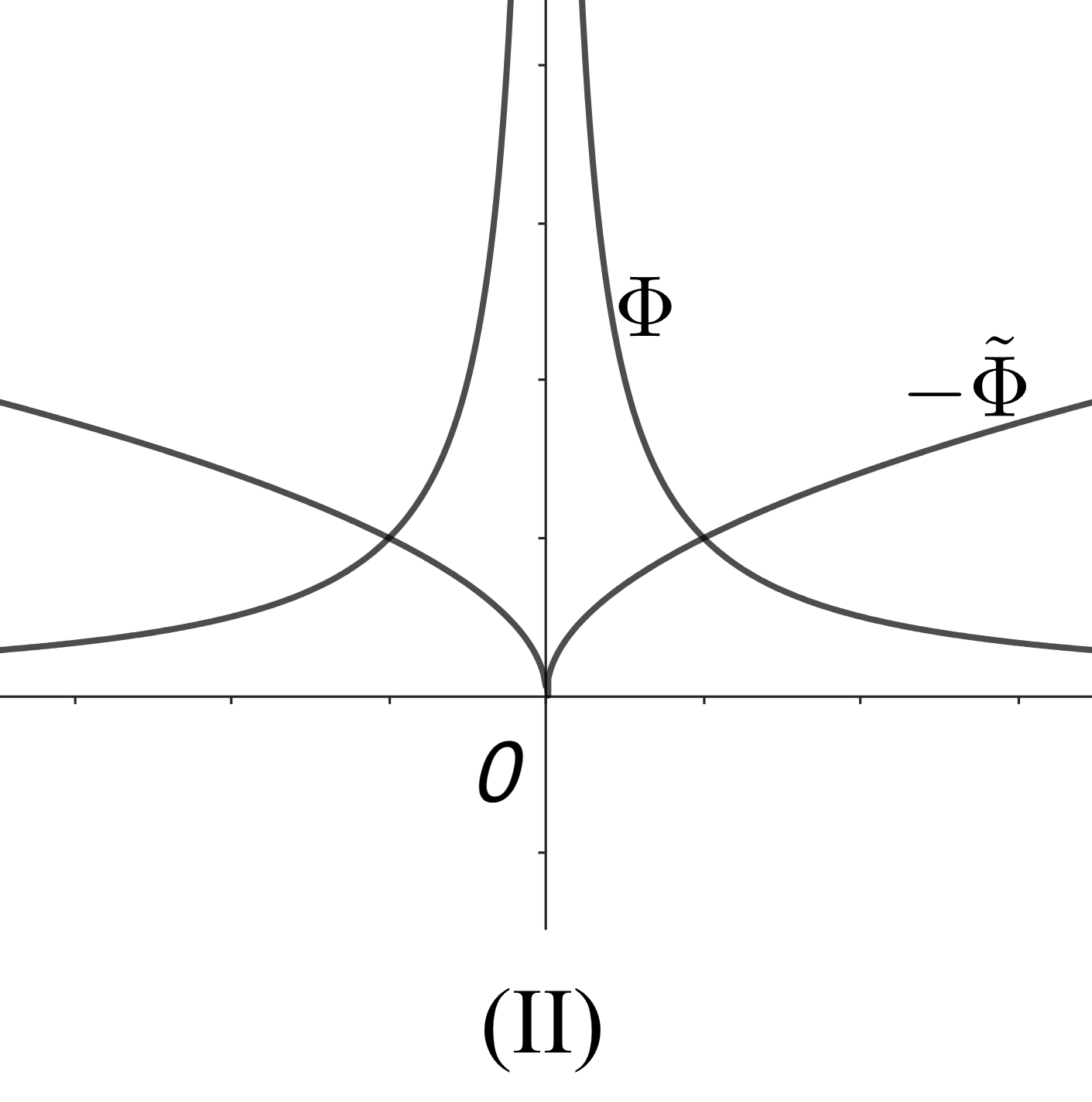}\\
  \includegraphics[width=60 mm]{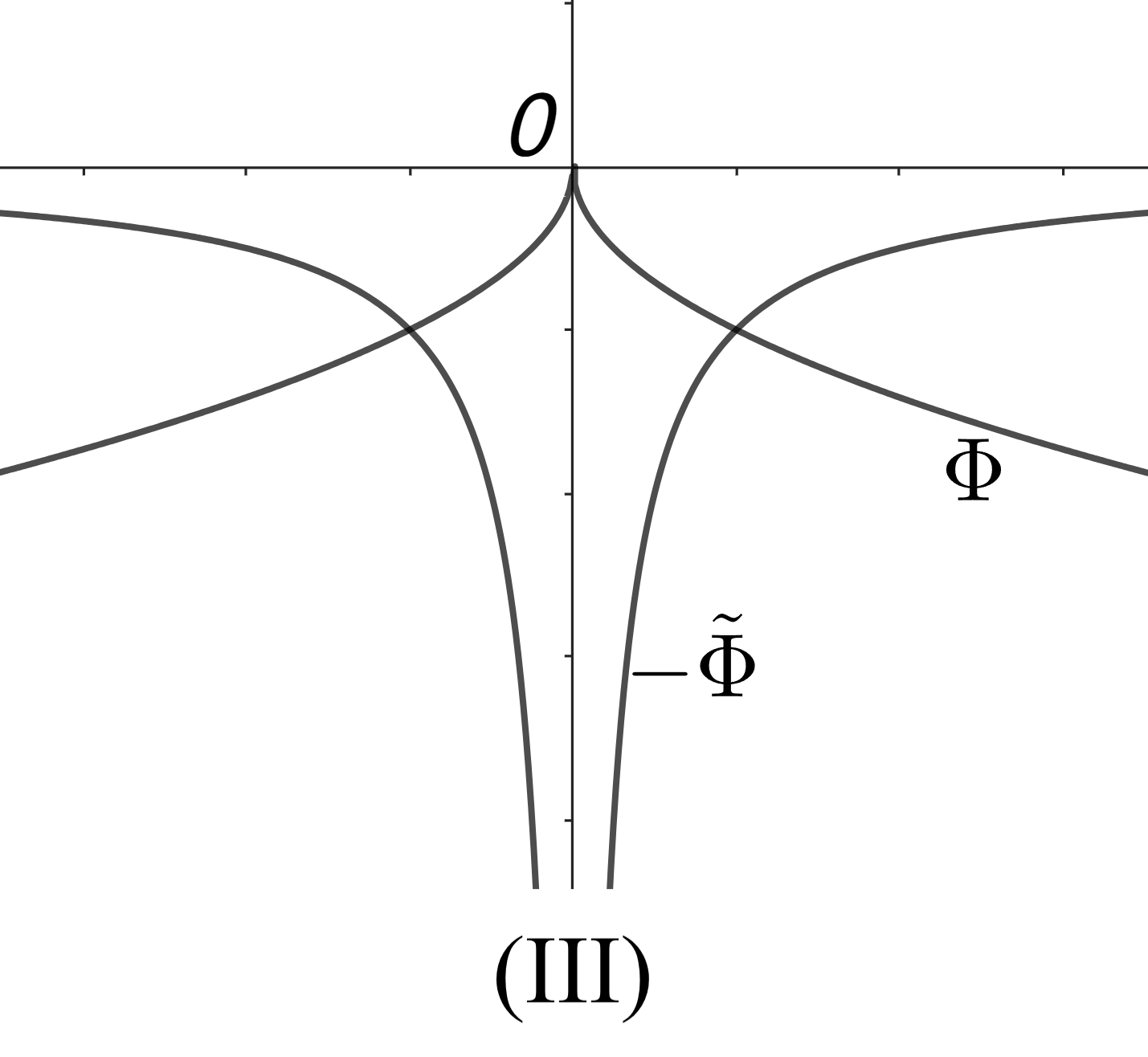}
  \includegraphics[width=60 mm]{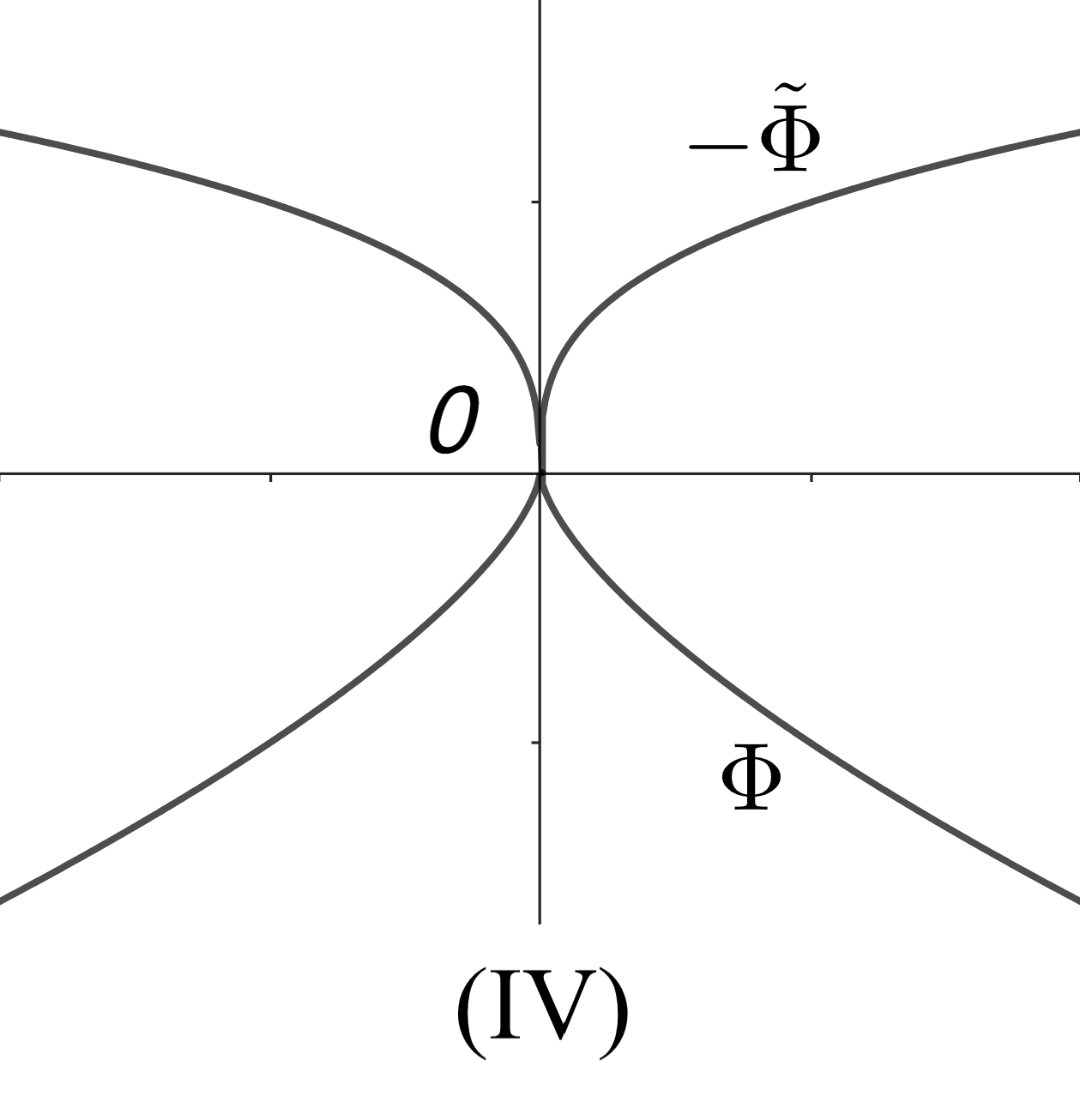}
  \caption{(I): $\alpha^*>0,\tilde{\alpha}^*>0$; (II): $\alpha^*>0,\tilde{\alpha}^*<0$;
  (III): $\alpha^*<0,\tilde{\alpha}^*>0$; (IV): $\alpha^*<0,\tilde{\alpha}^*<0$.}
\end{figure}

For general operators, there are four possible types of fundamental solution pairs $(\Phi, -\tilde{\Phi})$ due to the signs of $\alpha^*$ and $\tilde{\alpha}^*$. Their behavior at the singular point (i.e. the origin $0$) and at infinity are illustrated in Figure \ref{f.1}. As pointed out in \cite[P. 760]{MR2663711}, we don't know any operator $F$ with $\alpha^*<0$ and $\tilde\alpha^*<0$. Note that a fundamental solution could tend to infinity as $x\to \infty$ with the speed strictly less than linear growth since $\alpha^*,\tilde{\alpha}^*>-1$, which is critical for proving the convergence of derivatives of solutions as $x\to \infty$ (see \Cref{th1.0} and \Cref{th1.1}).

We use the following notation adopted from \cite{MR2663711}:
\begin{definition}\label{de1.2}
Let $0<a<+\infty$ and we write
\begin{equation*}
  \begin{aligned}
&(i)~ u\sim v~~\mbox{ if}~ u,v\to 0 ~\mbox{ and }~\frac{u}{v}\to a~~\mbox{as}~x\to \infty;\\
&(ii)~ u\approx v~~\mbox{ if}~ u,v\to \infty~\mbox{ and }~\frac{u}{v}\to a ~~\mbox{as}~x\to \infty.\\
  \end{aligned}
\end{equation*}
\end{definition}

\begin{remark}\label{re1.2}
In \cite{MR2663711}, (ii) is characterized as:
\begin{equation*}
u\approx v~~\mbox{ if}~ u,v\to \infty ~~\mbox{as}~x\to \infty~\mbox{ and }~av-C\leq u\leq av+C~~\mbox{ in }~B_1^c,
\end{equation*}
where $C$ is a constant. In fact, if $u,v$ are two solutions of \cref{e.F1} with $F$ being positively homogenous of degree $1$, these two definitions are equivalent. Indeed, since $u,v\to \infty$, the latter imply the former obviously. For the inverse, we assume that $u,v\to +\infty$ without loss of generality. Since $u/v\to a$, for any $0<\varepsilon<a/4$, there exists $R_0\geq 1$ such that for any $|x|\geq R_0$,
\begin{equation*}
(a-\varepsilon)v(x)\leq u(x)\leq (a+\varepsilon)v(x).
\end{equation*}
Since $u,v\to +\infty$, there exists $R_1\geq R_0$ such that for any $|x|\geq R_1$,
\begin{equation*}
(a-2\varepsilon)\left(v(x)-\max_{\partial B_1} v\right)+\min_{\partial B_1} u\leq u(x)
\leq (a+2\varepsilon)\left(v(x)-\min_{\partial B_1} v\right)+\max_{\partial B_1} u.
\end{equation*}
For any $R\geq R_1$, by the comparison principle (see \cite[Theorem 3.3]{MR1118699} or \cite[Corollary 3.7 and Theorem 5.3]{MR1351007},
\begin{equation*}
(a-2\varepsilon)\left(v-\max_{\partial B_1} v\right)+\min_{\partial B_1} u\leq u
\leq (a+2\varepsilon)\left(v-\min_{\partial B_1} v\right)+\max_{\partial B_1} u~~\mbox{ in}~B_R\backslash B_1.
\end{equation*}
By letting $R\to \infty$ and then $\varepsilon\to 0$,
\begin{equation*}
av-a\max_{\partial B_1} v+\min_{\partial B_1} u\leq u
\leq av-a\min_{\partial B_1} v+\max_{\partial B_1} u~~\mbox{ in}~B^c_1.
\end{equation*}
\end{remark}
~\\

A basic fact about the asymptotic behavior is that once a solution $u$ of \cref{e.F1} is bounded above or below, $\lim_{x\to \infty} u(x)$ exists (may be $+\infty$ or $-\infty$) and we denote it by $u_\infty$ (see \Cref{le2.1}). Armstrong, Sirakov and Smart \cite[Theorem 1.10]{MR2663711} have given the following complete classification of the asymptotic behavior at infinity if $u$ is bounded below or above. Roughly speaking, $u$ is comparable to one of the fundamental solutions $\Phi$ or $-\tilde{\Phi}$ at infinity (up to modulo $u_{\infty}$), or $u_{\infty}$ is attained on each sphere $\partial B_R$ ($R\geq 1$).
\begin{lemma}\label{pr2.3}
Let $u$ be a viscosity solution of \cref{e.F1} and be bounded above or below. Suppose that $F$ is positively homogeneous of degree $1$. Then precisely one of the following five alternatives holds:
\begin{flalign*}
\begin{split}
&(i)~\min_{\partial B_R}u\leq u_{\infty}\leq \max_{\partial B_R}u,~\forall ~R\geq1;\\
&(ii)~\alpha^*>0,~u-u_{\infty}\sim \Phi;\\
&(iii)~\tilde{\alpha}^*>0,~u-u_{\infty}\sim -\tilde\Phi;\\
&(iv)~\alpha^*\leq 0,~u\approx \Phi;\\
&(v)~\tilde{\alpha}^*\leq 0,~u\approx -\tilde\Phi.
\end{split}&
\end{flalign*}

In particular, if $\alpha^*>0$ and $\tilde{\alpha}^*>0$, $u_{\infty}$ is finite and for any $x\in B_1^c$,
\begin{equation*}
|u(x)-u_{\infty}|\leq C\|u-u_{\infty}\|_{L^{\infty}(\partial B_1)}\max(\Phi(x),\tilde{\Phi}(x))\leq C\|u-u_{\infty}\|_{L^{\infty}(\partial B_1)}|x|^{-\min(\alpha^*,\tilde\alpha^*)},
\end{equation*}
where $C$ is universal.
\end{lemma}

\begin{remark}\label{re1.4}
This classification is much more complicated than that of Lapalce equation. For the Laplace equation, $\Phi=\tilde{\Phi}=|x|^{2-n}$ ($n\geq 3$), which corresponds to (I) in Figure \ref{f.1}. However, for general nonlinear operators, there are four possible types of fundamental solution pairs. Hence, the classification is complicated.

For the reader's convenience, we give the proof of this lemma (see the end of \Cref{S2}).
\end{remark}

\begin{remark}\label{re2.1}
If we consider the Laplace equation and $n\geq 3$, then $\alpha^*=\tilde{\alpha}^*=n-2>0$. Then $u_{\infty}$ is finite. Consider the Kelvin transformation:
\begin{equation*}
v(x):=|x|^{2-n}\left(u\left(\frac{x}{|x|^2}\right)-u_{\infty}\right).
\end{equation*}
Hence, $v$ is harmonic in $B_1\backslash \{0\}$ and $v=o(|x|^{2-n})$ as $x\to 0$. Thus, $0$ is a removable singularity. Therefore, for any $k\geq 0$, there exists constants $a_{\sigma}$ such that
\begin{equation*}
\left|v(x)-\sum_{|\sigma|=0}^{k} a_{\sigma}x^{\sigma}\right|\leq C\|v\|_{L^{\infty}(\partial B_1)}|x|^{k+1},~\forall ~x\in B_1,
\end{equation*}
where $\sigma$ is the multi-index and $C$ depends only on $n,\lambda,\Lambda$ and $k$. By transforming to $u$, we have
\begin{equation*}
\left|u(x)-u_{\infty}-|x|^{2-n}\sum_{|\sigma|=0}^{k} a_{\sigma}x^{\sigma}|x|^{-2|\sigma|}\right|
\leq C\|u-u_{\infty}\|_{L^{\infty}(\partial B_1)}|x|^{-n-k+1},~\forall ~x\in B_1^c.
\end{equation*}
\end{remark}
~\\

Now, we state our main results.
\begin{theorem}\label{th1.0}
Let $u\in C(B_1^c)$ be a viscosity solution of \cref{e.F1}. Suppose that $F(0)=0$ and
\begin{equation}\label{e1.2-2}
|u(x)|\leq |x|~~\mbox{ in}~~B_1^c.
\end{equation}
Then there exists a linear polynomial $P$ such that~\\
(i) if $-1<\alpha^+<0$, i.e., $\Lambda/\lambda> n-1$,
\begin{equation}\label{e2.1-2}
\begin{aligned}
&|u(x)-P(x)|\leq C|E^+(x)|=C|x|^{1-(n-1)\lambda/\Lambda},~\forall ~x\in B_2^c,\\
&|Du(x)-DP|\leq C|x|^{-1}|E^+(x)|=C|x|^{-(n-1)\lambda/\Lambda},~\forall ~x\in B_4^c;
\end{aligned}
\end{equation}
(ii) if $\alpha^+=0$, i.e., $\Lambda/\lambda= n-1$,
\begin{equation}\label{e2.1-3}
\begin{aligned}
&|u(x)-P(x)|\leq C\ln |x|,~\forall ~x\in B_2^c,\\
&|Du(x)-DP|\leq C|x|^{-1}\ln |x|,~\forall ~x\in B_4^c;
\end{aligned}
\end{equation}
(iii) if $\alpha^+>0$, i.e., $\Lambda/\lambda< n-1$,
\begin{equation}\label{e4.1}
\begin{aligned}
&|u(x)-P(x)|\leq C\|u-P\|_{L^{\infty}(\partial B_1)}|x|^{1-(n-1)\lambda/\Lambda},~\forall ~x\in B_2^c,\\
&|Du(x)-DP|\leq C\|u-P\|_{L^{\infty}(\partial B_1)}|x|^{-(n-1)\lambda/\Lambda},~\forall ~x\in B_4^c,
\end{aligned}
\end{equation}
where $C$ is universal.
\end{theorem}

\begin{remark}\label{re1.9}
Since we always have $\alpha^+>-1$, $Du(x)\to DP$ as $x\to \infty$. If $\alpha^+>0$, i.e., $\Lambda/\lambda< n-1$, we also have $u(x)-P(x)\to 0$ as $x\to \infty$. That is, the asymptotic behavior of $u$ at infinity is like a linear polynomial. However, if $\alpha^+\leq 0$, we only have a control of $u-P$.
\end{remark}

\begin{remark}\label{re1.5}
Note that $F(0)=0$ is necessary. For instance, we can't infer that $u$ is a linear polynomial if
\begin{equation*}
F(D^2u)=\Delta u-1=0~~\mbox{ in}~\mathbb{R}^n
\end{equation*}
and
\begin{equation*}
|u(x)|\leq |x|+1~~\mbox{ in}~~\mathbb{R}^n.
\end{equation*}
\end{remark}

\begin{remark}\label{re1.8}
A general exterior problem $F(D^2u)=0$ with $|u(x)|\leq K(|x|+1)$ in $\Omega^c$ (a compact set $\Omega$ and a constant $K$) can be transformed to \cref{e.F1} with \cref{e1.2-2} easily. Hence, we treat the latter for simplicity in this paper.
\end{remark}
~\\

With the aid of \Cref{pr2.3}, we have a complete classification of the asymptotic behavior at infinity if the solution grows at most linearly.
\begin{corollary}\label{th1.0-2}
Let $u\in C(B_1^c)$ be a viscosity solution of \cref{e.F1}. Suppose that $F$ is positively homogeneous of degree $1$ and
\begin{equation*}
|u(x)|\leq |x|~~\mbox{ in}~~B_1^c.
\end{equation*}
Then there exists a linear polynomial $P$ such that precisely one of the following five alternatives holds:
\begin{flalign*}
\begin{split}
&(i)~\min_{\partial B_R}(u-P)\leq 0\leq \max_{\partial B_R}(u-P),~\forall ~R\geq 1;\hfill\\
&(ii)~\alpha^*>0,~u-P\sim \Phi;\\
&(iii)~\tilde{\alpha}^*>0,~u-P\sim -\tilde\Phi;\\
&(iv)~\alpha^*\leq 0,~u-P\approx \Phi;\\
&(v)~\tilde{\alpha}^*\leq 0,~u-P\approx -\tilde\Phi.
\end{split}&
\end{flalign*}

In particular, if $\alpha^*>0$ and $\tilde{\alpha}^*>0$, then
\begin{equation}\label{e4.2}
|u(x)-P(x)|\leq C\|u-P\|_{L^{\infty}(\partial B_1)}|x|^{-\min(\alpha^*,\tilde\alpha^*)},
\end{equation}
where $C$ is universal.
\end{corollary}

Next, we consider the asymptotic behavior of solutions with quadratic growth.
\begin{theorem}\label{th1.1}
Let $u\in C(B_1^c)$ be a viscosity solution of
\begin{equation}\label{e.F3}
F(D^2u)=A~~\mbox{ in}~~\bar B_1^c,
\end{equation}
where $A\in \mathbb{R}$. Suppose that $F$ is a convex function and
\begin{equation}\label{e1.2}
|u(x)|\leq |x|^2~~\mbox{ in}~~B_1^c.
\end{equation}
Then there exists a quadratic polynomial $P$ with $F(D^2P)=A$ such that~\\
(i) if $-1<\alpha^+<0$, i.e., $\Lambda/\lambda> n-1$,
\begin{equation}\label{e2.1-2-2}
\begin{aligned}
&|u(x)-P(x)|\leq C|E^+(x)|=C|x|^{1-(n-1)\lambda/\Lambda},~\forall ~x\in B_2^c,\\
&|Du(x)-DP|\leq C|x|^{-1}|E^+(x)|=C|x|^{-(n-1)\lambda/\Lambda},~\forall ~x\in B_4^c,\\
&|D^2u(x)-D^2P|\leq C|x|^{-2}|E^+(x)|=C|x|^{-1-(n-1)\lambda/\Lambda},~\forall ~x\in B_4^c,
\end{aligned}
\end{equation}
where $C$ depends only on $n,\lambda,\Lambda,|A|$ and $|F(0)|$;\\
(ii) if $\alpha^+=0$, i.e., $\Lambda/\lambda= n-1$,
\begin{equation}\label{e2.1-3-2}
\begin{aligned}
&|u(x)-P(x)|\leq C\ln |x|,~\forall ~x\in B_2^c,\\
&|Du(x)-DP|\leq C|x|^{-1}\ln |x|,~\forall ~x\in B_4^c,\\
&|D^2u(x)-D^2P|\leq C|x|^{-2}\ln |x|,~\forall ~x\in B_4^c,
\end{aligned}
\end{equation}
where $C$ depends only on $n,\lambda,\Lambda,|A|$ and $|F(0)|$;\\
(iii) if $\alpha^+>0$, i.e., $\Lambda/\lambda< n-1$,
\begin{equation}\label{e4.1-2}
\begin{aligned}
&|u(x)-P(x)|\leq C\|u-P\|_{L^{\infty}(\partial B_1)}|x|^{1-(n-1)\lambda/\Lambda},~\forall ~x\in B_2^c,\\
&|Du(x)-DP|\leq C\|u-P\|_{L^{\infty}(\partial B_1)}|x|^{-(n-1)\lambda/\Lambda},~\forall ~x\in B_4^c,\\
&|D^2u(x)-D^2P|\leq C\|u-P\|_{L^{\infty}(\partial B_1)}|x|^{-1-(n-1)\lambda/\Lambda},~\forall ~x\in B_4^c,
\end{aligned}
\end{equation}
where $C$ is universal.
\end{theorem}

\begin{remark}\label{re1.6}
The condition ``$F$ is convex'' is used to guarantee the interior $C^{2,\alpha}$ regularity (see \cite[Theorem 6.6]{MR1351007}). It can be removed in dimension two (see \cite[Theorem 4.9]{MR4560756}). Even in higher dimension, this condition can be relaxed to some extent (see \cite{MR1995493} and \cite{MR1793687}).
\end{remark}

\begin{remark}\label{re1.7}
Note that we don't have a complete classification as \Cref{th1.0-2}. The reason is the following. Since $P$ is a quadratic polynomial, $v:=u-P$ is a solution of
\begin{equation*}
\tilde{F}(D^2v):=F(D^2v+D^2P)=A~~\mbox{ in}~\bar{B}_1^c.
\end{equation*}
Hence, $\tilde{F}$ is not positively homogeneous of degree $1$ in general even if $F$ is.
\end{remark}

\begin{remark}\label{re1.10}
Li, Li and Yuan proved a similar result (see \cite[Theorem 2.1]{MR4038557}) with a different approach and the smoothness of $F$ was assumed.
\end{remark}
~\\

\section{Preliminaries}\label{S2}
In this section, we prepare some preliminaries. Throughout this section, we always assume that $F$ is positively homogeneous of degree $1$. Most of contents can be founded in \cite{MR2663711}. For completeness and the reader's convenience, we give the detailed proofs. The idea are mainly inspired by \cite{MR2663711} and \cite{MR4038557}.

The first fact is that once a solution is bounded from one side, the limit at infinity exists. Precisely, we have
\begin{lemma}\label{le2.1}
Let $u\in S(\lambda,\Lambda,0)$ and be bounded below or above. Then $u_{\infty}:=\lim_{x\rightarrow \infty} u(x)$ exists (may be $+\infty$ or $-\infty$).
\end{lemma}
\proof We only prove the case that $u$ is bounded below. Indeed, if $u$ is bounded above, we can consider $-u\in S(\lambda,\Lambda,0)$. Without loss of generality, we assume that $u\geq 0$.

Let $a=\liminf_{x\rightarrow \infty} u(x)$. If $a=+\infty$, $\lim_{x\rightarrow \infty} u=+\infty$ clearly. If $a<+\infty$, for any $\varepsilon>0$, there exist two sequences of $R_i\rightarrow \infty$ increasingly and $x_i\in \partial B_{R_i}$ such that $R_2\geq 2R_1$,
\begin{equation*}
u-a+\varepsilon\geq 0~~\mbox{ in } \bar B^c_{R_1} ~~\mbox{ and }~~u(x_i)-a-\varepsilon\leq 0,~\forall ~i\geq 2.
\end{equation*}
By applying the Harnack inequality (see \cite[Theorem 4.3]{MR1351007}) on each $\partial B_{R_i}$ ($i\geq 2$),
\begin{equation}\label{e2.4}
u(x)-a+\varepsilon\leq C(u(x_i)-a+\varepsilon)\leq C\varepsilon,~\forall ~x\in \partial B_{R_i},
\end{equation}
where $C>0$ is universal.

Next, for any $x\in \bar B^c_{R_2}$, there exists $i\geq 2$ such that $x\in B_{R_{i+1}}\backslash B_{R_{i}}$. By the comparison principle and \cref{e2.4},
\begin{equation*}
u(x)\leq a+C\varepsilon.
\end{equation*}
Therefore, we have proved that for any $\varepsilon>0$, there exists $R_2>0$ such that
\begin{equation*}
-\varepsilon \leq u(x)-a\leq 2C\varepsilon,~\forall ~x\in B^c_{R_2}.
\end{equation*}
That is, $\lim_{x\rightarrow \infty} u=a$.~\qed
~\\

The next two lemmas demonstrate how the behavior of a solution can be influenced by another solution. This is why the fundamental solutions play a key role in the study of asymptotic behavior.
\begin{lemma}\label{le2.2}
(i) If $v$ is a supersolution and $w$ is a subsolution of \cref{e.F1} with
\begin{equation*}
\lim_{x\rightarrow \infty} v=+\infty,~\limsup_{x\rightarrow \infty}\frac{w(x)}{v(x)}\leq 0,
\end{equation*}
then $w$ is bounded above and
\begin{equation}\label{e3.1}
w(x)\leq \max_{\partial B_R} w,~\forall ~x\in B_R^c,~\forall ~R\geq 1.
\end{equation}

(ii) Similarly, if $v$ is a subsolution and $w$ is a supersolution of \cref{e.F1} with
\begin{equation*}
\lim_{x\rightarrow \infty} v=-\infty,~\limsup_{x\rightarrow \infty}\frac{w(x)}{v(x)}\leq 0,
\end{equation*}
then $w$ is bounded below and
\begin{equation}\label{e3.2}
w(x)\geq \min_{\partial B_R} w,~\forall ~x\in B_R^c,~\forall ~R\geq 1.
\end{equation}
\end{lemma}
\proof We only prove the first case since the second case can be proved in a similar way. Without loss of generality, we assume $v> 0$. By the assumption, for any $\varepsilon>0,R\geq1$, there exists $R_0>1$ such that for any $R_1\geq R_0$,
\begin{equation*}
w\leq \varepsilon v/2~~\mbox{ and }~~\max_{\partial B_R}w\geq- \varepsilon v/2~~\mbox{ on}~~\partial B_{R_1}.
\end{equation*}
Hence,
\begin{equation*}
w\leq \max_{\partial B_R}w+\varepsilon v~~\mbox{ on }~~\partial B_{R_1}.
\end{equation*}
By the comparison principle (note that $v>0$),
\begin{equation*}
w\leq \max_{\partial B_R} w+\varepsilon v~~\mbox{ in }~~B_{R_1}\backslash B_R.
\end{equation*}
Let $R_1\rightarrow \infty$ first and then $\varepsilon\rightarrow 0$, we have
\begin{equation*}
  w\leq \max_{\partial B_R}w~~\mbox{ in }~~B_R^c.
\end{equation*}
~\qed~\\

\begin{lemma}\label{le2.2-2}
Suppose that
\begin{equation}\label{e3.3}
\lim_{x\rightarrow \infty} v(x)=\lim_{x\rightarrow \infty} w(x)=0.
\end{equation}
(i) If $v$ is a positive supersolution and $w$ is a subsolution of \cref{e.F1}, then there exists $a\geq 0$ such that
\begin{equation*}
w\leq av~~\mbox{ in}~~B_1^c.
\end{equation*}
(ii) If $v$ is a subsolution and $w$ is a supersolution with $w> 0$ on $\partial B_R$ for some $R\geq 1$, then there exists $a> 0$ such that
\begin{equation*}
w\geq a v~\mbox{ in}~~B_R^c.
\end{equation*}
(iii) If $v$ is a negative subsolution and $w$ is a supersolution, then there exists $a\geq 0$ such that
\begin{equation*}
w\geq a v~~\mbox{ in}~~B_1^c,
\end{equation*}
(iv) If $v$ is a supersolution and $w$ is a subsolution with $w<0$ on $\partial B_R$ for some $R\geq 1$, then there exists $a> 0$ such that
\begin{equation*}
w\leq av~\mbox{ in}~~B_R^c.
\end{equation*}
\end{lemma}
\proof We only prove (i) and (ii). For (i), since $v$ is positive, there exists $a\geq 0$ such that $w\leq av$ on $\partial B_1$. Note that $\lim_{x\rightarrow \infty} v=\lim_{x\rightarrow \infty} w=0$. By the comparison principle, for any $\varepsilon>0$, there exists $R_0>0$ such that for any $R\geq R_0$,
\begin{equation*}
w\leq av+\varepsilon~~\mbox{ in}~~B_R\backslash B_{1}.
\end{equation*}
Let $R\rightarrow \infty$ first and $\varepsilon\rightarrow 0$ next. Then
\begin{equation*}
w\leq av~~\mbox{ in}~~B_1^c.
\end{equation*}
That is, (i) holds.

Next we prove (ii) in a similar way. Since $w> 0$ on $\partial B_R$, there exists $a> 0$ such that $w\geq av$ on $\partial B_R$. By the comparison principle as above, for any $\varepsilon>0$, there exists $R_0>0$ such that for any $R_1\geq R_0$,
\begin{equation*}
w\geq a v-\varepsilon~~\mbox{ in}~~B_{R_1}\backslash B_{R}.
\end{equation*}
Let $R_1\rightarrow \infty$ first and $\varepsilon\rightarrow 0$ next. Then
\begin{equation*}
w\geq av~~\mbox{ in}~~B_R^c.
\end{equation*}
~\qed
~\\

Next, we show that the limit of the ratio of two solutions exists if the lower limit exists, which is a simple consequence of the Harnack inequality and the comparison principle.
\begin{lemma}\label{le2.8}
Suppose that $v,w$ are solutions of \cref{e.F1} with $w> 0$ (or $w<0$) and
\begin{equation*}
0<a:=\liminf_{x\rightarrow \infty} \frac{v(x)}{w(x)}<+\infty.
\end{equation*}
Then
\begin{equation}\label{e2.13}
\lim_{x\rightarrow \infty}\frac{v(x)}{w(x)}=a.
\end{equation}

If $v>0$ or $v<0$ in addition and $a=0$, then
 \begin{equation}\label{e2.14}
\lim_{x\rightarrow \infty}\frac{v(x)}{w(x)}=0.
\end{equation}
\end{lemma}
\proof We only consider the case $w>0$. By the definition of $a$, for any $0<\varepsilon<a/2$, there exist sequences of $R_i\rightarrow \infty$ increasingly and $x_i\in \partial B_{R_i}$ such that $R_2\geq 2R_1$,
\begin{equation*}
v-(a-\varepsilon)w\geq 0~~\mbox{ in } B^c_{R_1} ~~\mbox{ and }~~
v(x_i)-(a+\varepsilon)w(x_i) \leq 0.
\end{equation*}
Since $F$ is positively homogeneous of degree $1$ and $a-\varepsilon>0$, $(a-\varepsilon)w$ is also a solution and then
\begin{equation*}
v-(a-\varepsilon)w\in S(\lambda,\Lambda,0)~~\mbox{ in}~~ B^c_{1}.
\end{equation*}
Then by applying the Harnack inequality to each $\partial B_{R_i}$ ($i\geq 2$),
\begin{equation*}
v(x)-(a-\varepsilon)w(x)
\leq C\left(v(x_i)-(a-\varepsilon)w(x_i) \right)
\leq C\varepsilon w(x_i)
\leq C\varepsilon w(x),~\forall ~x\in \partial B_{R_i}.
\end{equation*}
That is,
\begin{equation*}
 v\leq (a+C\varepsilon)w~~\mbox{ on}~~\partial B_{R_i},~\forall ~i\geq 2.
\end{equation*}
From the comparison principle,
\begin{equation*}
v\leq (a+C\varepsilon)w~~\mbox{ in}~~B_{R_i}\backslash B_{R_2}.
\end{equation*}
Hence,
\begin{equation*}
(a-\varepsilon)w\leq v\leq (a+C\varepsilon)w~~\mbox{ in}~~B_{R_i}\backslash B_{R_2}.
\end{equation*}
By taking $i\rightarrow \infty$,
\begin{equation*}
-\varepsilon\leq \frac{v}{w}-a\leq C\varepsilon~~\mbox{ in}~~B^c_{R_2}.
\end{equation*}
Hence, we have proved that for any $\varepsilon>0$, there exists $R_2>0$ such that
\begin{equation*}
\left|\frac{v(x)}{w(x)}-a\right|\leq C\varepsilon,~~\forall ~ x\in B^c_{R_2}.
\end{equation*}
That is, \cref{e2.13} holds.

Finally, we consider $v,w>0$ and $a=0$ (other cases can be treated similarly). By the Harnack inequality,
\begin{equation*}
v(x)\leq Cv(x_i)\leq C\varepsilon w(x_i)
\leq C\varepsilon w(x),~\forall ~x\in \partial B_{R_i},~\forall ~i\geq 2,
\end{equation*}
where $\varepsilon>0$, $R_i$ and $x_i$ are as above. From the comparison principle,
\begin{equation*}
v\leq C\varepsilon w~~\mbox{ in}~~B_{R_i}\backslash B_{R_2},~\forall ~i\geq 3.
\end{equation*}
By letting $i\to \infty$, we infer that
\begin{equation*}
\limsup_{x\rightarrow \infty} \frac{v(x)}{w(x)}\leq C\varepsilon.
\end{equation*}
Since $\varepsilon$ is arbitrary, \cref{e2.14} holds.~\qed ~\\

Now, we study the behavior for functions belonging to the Pucci's class, which will be used in next section.
\begin{lemma}\label{le2.3}
Let $u\in S(\lambda,\Lambda,0)$ in $\bar{B}_1^c$ and $u\geq 0$. Then
\begin{flalign*}
  \begin{split}
&(i) ~\mbox{if}~ \alpha^+>0,~|u(x)-u_{\infty}|\leq \|u-u_{\infty}\|_{L^{\infty}(\partial B_1)}E^+(x)~~\mbox{ in}~~B_2^c;\\
&(ii) ~\mbox{if}~\alpha^+\leq 0,~ 0\leq u(x)\leq -C\|u\|_{L^{\infty}(B_4\backslash B_1)}E^+(x)~~\mbox{ in}~~B_2^c,
  \end{split}&
\end{flalign*}
where $C$ is universal.


Similarly, if $u\leq 0$, then
\begin{flalign*}
  \begin{split}
&(i) ~\mbox{if}~ \alpha^+>0,~|u(x)-u_{\infty}|\leq \|u-u_{\infty}\|_{L^{\infty}(\partial B_1)}E^+(x)~~\mbox{ in}~~B_2^c;\\
&(ii) ~\mbox{if}~\alpha^+\leq 0,~ C\|u\|_{L^{\infty}(B_4\backslash B_1)}E^+(x)\leq u(x)\leq 0~~\mbox{ in}~~B_2^c.
  \end{split}&
\end{flalign*}
\end{lemma}
\proof We only give the proof for $u\geq 0$. We first prove (i). By \Cref{le2.1}, $u_{\infty}$ is well defined. Note that $u_\infty<+\infty$. Otherwise, by \Cref{le2.2} (i) with $F=\mathcal{M}^-$, $v=u$ and $w=e^-=-E^+$, we have
\begin{equation*}
-E^+(x)\leq \max_{\partial B_1} (-E^+)=-1,~\forall ~x\in B_1^c,
\end{equation*}
which is impossible. Next, by applying \Cref{le2.2-2} (i) (with $F=\mathcal{M}^+,v=E^+,w=u-u_\infty$) and (iii) (with $F=\mathcal{M}^-, v=e^-=-E^+,w=u-u_\infty$) respectively, the conclusion follows.

Next, we prove (ii). Since $\alpha^+\leq 0$, $e^-(x)=-E^+(x)\rightarrow +\infty$ increasingly as $|x|\to \infty$. Let
\begin{equation*}
a=\sup_{x\in B_2^c} \frac{u(x)}{e^-(x)}.
\end{equation*}
Clearly, $a>0$. If $a=+\infty$, for any $K>0$, there exist sequences of $R_i\rightarrow \infty$ increasingly and $x_i\in \partial B_{R_i}$ such that $R_2\geq 2R_1$ and
\begin{equation*}
u(x_i)\geq Ke^-(x_i).
\end{equation*}
By applying the Harnack inequality to $u$ on $\partial B_{R_i}$,
\begin{equation*}
u(x)\geq cu(x_i)\geq cKe^-(x_i)=cKe^-(x),~\forall ~i\geq 2,~\forall ~x\in \partial B_{R_i},
\end{equation*}
where $c>0$ is universal.

Next, by the comparison principle, for any $i\geq 3$,
\begin{equation*}
u \geq cKe^-~~\mbox{ in}~ B_{R_i}\backslash B_{R_2}.
\end{equation*}
Let $i\to \infty$ and we have
\begin{equation*}
u\geq cKe^-~~\mbox{ in}~ B_{R_2}^c.
\end{equation*}
Therefore, we have proved that for any $K>0$, there exists $R_2>0$ such that
\begin{equation*}
u(x)\geq cKe^-(x),~\forall ~x\in B_{R_2}^c.
\end{equation*}
That is, $\lim_{x\rightarrow \infty} u(x)/e^-(x)=+\infty$. In another word, $\lim_{x\rightarrow \infty} e^-(x)/u(x)=0$. By \Cref{le2.2} (i) with $F=\mathcal{M}^-$, $v=u$ and $w=e^-$, we arrive at a contradiction. Hence, $a<+\infty$.

By the definition of $a$, there exists $R_0\geq 2$ and $x_0\in \partial B_{R_0}$ such that $u(x_0)\geq ae^-(x_0)/2$. If $R_0\leq 4$,
\begin{equation}\label{e5.2}
a\leq \frac{2u(x_0)}{e^-(x_0)}\leq C \|u\|_{L^{\infty}(B_4\backslash B_1)}.
\end{equation}
If $R_0\geq 4$, by the Harnack inequality,
\begin{equation*}
  u\geq cae^-~~\mbox{ on}~\partial B_{R_0},
\end{equation*}
where $c>0$ is universal. Next, the comparison principle implies
\begin{equation*}
ca+u\geq cae^-~~\mbox{ in}~B_{R_0}\backslash B_1.
\end{equation*}
Hence, by applying above inequality on $\partial B_4$ and noting $e^-\geq C_0>1$ on $\partial B_4$,
\begin{equation*}
ca+\|u\|_{L^{\infty}(\partial B_4)}\geq C_0ca.
\end{equation*}
Thus,
\begin{equation}\label{e5.3}
a\leq C\|u\|_{L^{\infty}(\partial B_4)}\leq C\|u\|_{L^{\infty}(B_4\backslash B_1)}.
\end{equation}
Finally, by combining \cref{e5.2} and \cref{e5.3}, we obtain (ii).~\qed~\\

%

At the end of this section, For the completeness and the reader's convenience, we prove the classification of the asymptotic behavior at infinity:~\\
\noindent\textbf{Proof of \Cref{pr2.3}.} We only prove the case that $u$ is bounded below. Without loss of generality, we assume $u\geq 0$.

By \Cref{le2.1}, $u_{\infty}:=\lim_{x\rightarrow \infty} u(x)$ exists. First, we consider the case $u_{\infty}=+\infty$. If $\tilde{\alpha}^*>0$,
\begin{equation}\label{e2.5}
-\tilde{\Phi}(x)\to 0~~\mbox{increasingly as}~ x\to \infty.
\end{equation}
By \Cref{le2.2} (i) with $v=u, w=-\tilde{\Phi}$,
\begin{equation*}
-\tilde{\Phi}(x)\leq \max_{\partial B_1}-\tilde{\Phi}=-1,~\forall ~x\in B_1^c,
\end{equation*}
which contradicts with \cref{e2.5}. Hence, $\tilde{\alpha}^*\leq 0$ and $\lim_{x\to \infty} -\tilde{\Phi} =+\infty$. Consider
\begin{equation*}
a:=\liminf_{x\rightarrow \infty}\frac{u(x)}{-\tilde{\Phi}(x)}.
\end{equation*}
If $a=0$ or $+\infty$, \Cref{le2.8} and \Cref{le2.2} (i) imply that one of $u$ and $-\tilde{\Phi}$ is bounded above, which is impossible. Hence, $0<a<+\infty$. Then by \Cref{le2.8} again,
\begin{equation*}
\lim_{x\rightarrow \infty}\frac{u(x)}{-\tilde{\Phi}(x)}=a.
\end{equation*}
That is, $u\approx -\tilde\Phi$.

If $u_{\infty}<+\infty$, we consider
\begin{equation*}
b:=\liminf_{x\to \infty} \frac{u(x)-u_{\infty}}{\Phi(x)},~c:=\liminf_{x\to \infty} \frac{u(x)-u_{\infty}}{-\tilde{\Phi}(x)}.
\end{equation*}
With the aid of \Cref{le2.2-2}, $b,c<+\infty$. If $b>0$, by \Cref{le2.8}, $u\sim \Phi$. Similarly, if $c>0$, $u\sim -\tilde\Phi$.

Finally, we consider: $b\leq 0$ and $c\leq 0$. We prove
\begin{equation}\label{e2.6}
\min_{\partial B_R}u\leq u_{\infty}\leq \max_{\partial B_R}u,~\forall ~R\geq 1
\end{equation}
for (i) $\alpha^*>0, \tilde{\alpha}^*>0$; (ii) $\alpha^*>0, \tilde{\alpha}^*\leq 0$; (iii) $\alpha^*\leq 0, \tilde{\alpha}^*>0$ and (iv) $\alpha^*\leq 0, \tilde{\alpha}^*\leq 0$ respectively.

For (i) $\alpha^*>0, \tilde{\alpha}^*>0$, we have $\Phi,-\tilde{\Phi}\to 0$ as $x\to \infty$. By \Cref{le2.2-2} (ii) and (iv), for any $R\geq 1$,
\begin{equation*}
\max_{\partial B_R}(u-u_{\infty})\geq 0~~\mbox{ and }~\min_{\partial B_R}(u-u_{\infty})\leq 0.
\end{equation*}
Otherwise, we have a contradiction with $b,c\leq 0$. Hence, \cref{e2.6} holds.

For (ii) $\alpha^*>0, \tilde{\alpha}^*\leq 0$, we have $\Phi\to 0$ and $-\tilde{\Phi}\to +\infty$ as $x\to \infty$. Since $b\leq 0$, with the aid of \Cref{le2.2-2} (ii), we obtain as above
\begin{equation*}
\min_{\partial B_R}(u-u_{\infty})\leq 0.
\end{equation*}
In addition, since $-\tilde{\Phi}\to +\infty$, by \Cref{le2.2} (i), we have
\begin{equation*}
u(x)\leq \max_{\partial B_R} u,~\forall ~x\in B_R^c,~\forall ~R\geq 1.
\end{equation*}
Let $x\to \infty$ and then $u_{\infty}\leq \max_{\partial B_R} u$. That is, \cref{e2.6} holds.

We can prove \cref{e2.6} for (iii) $\alpha^*\leq 0, \tilde{\alpha}^*>0$ and (iv) $\alpha^*\leq 0, \tilde{\alpha}^*\leq 0$ in a similar way and we omit the proof.~\qed~\\

\section{Asymptotic behavior at infinity}
In this section, we prove the asymptotic behavior at infinity for solutions with linear or quadratic growth. As pointed out in the introduction, this idea is highly motivated by \cite{MR4201294}. Since the proof uses the Liouville theorems in the whole space, let's recall them first.
\begin{lemma}\label{le3.1}
Let $u\in C(\mathbb{R}^n)$ be a viscosity solution of
\begin{equation*}
F(D^2u)=0~~\mbox{ in}~~\mathbb{R}^n.
\end{equation*}
Suppose that $F(0)=0$ and for some constant $K>0$,
\begin{equation*}
|u(x)|\leq K(|x|+1)~~\mbox{ in}~~\mathbb{R}^n.
\end{equation*}
Then $u$ is a linear polynomial.
\end{lemma}
\proof Since $F$ is uniformly elliptic, by the interior $C^{1,\alpha}$ regularity (see \cite[Corollary 5.7]{MR1351007}), there exists a linear polynomial $P$ such that for any $R>0$,
\begin{equation*}
|u(x)-P(x)|\leq C\frac{|x|^{1+\alpha}}{R^{1+\alpha}}\left(\|u\|_{L^{\infty}(B_{R})}+R^2F(0)\right),
~\forall ~x\in B_{R/2},
\end{equation*}
where $0<\alpha<1$ and $C>0$ are universal constants. Since $F(0)=0$ and $\|u\|_{L^{\infty}(B_{R})}\leq K(R+1)$, for any fixed $x\in \mathbb{R}^n$, by letting $R\to \infty$, we have
\begin{equation*}
  u(x)=P(x).
\end{equation*}
That is, $u$ is a linear polynomial.~\qed~\\

\begin{lemma}\label{le3.2}
Let $u\in C(\mathbb{R}^n)$ be a viscosity solution of
\begin{equation*}
F(D^2u)=A~~\mbox{ in}~~\mathbb{R}^n.
\end{equation*}
Suppose that $F$ is a convex function and for some constant $K>0$,
\begin{equation*}
|u(x)|\leq K(|x|^2+1)~~\mbox{ in}~~\mathbb{R}^n.
\end{equation*}
Then $u$ is a quadratic polynomial.
\end{lemma}
\proof Since $F$ is uniformly elliptic and convex, by the interior $C^{2,\alpha}$ regularity (see \cite[Theorem 6.6]{MR1351007}), there exists a quadratic polynomial $P$ such that for any $R>0$,
\begin{equation*}
|u(x)-P(x)|\leq C\frac{|x|^{2+\alpha}}{R^{2+\alpha}}\left(\|u\|_{L^{\infty}(B_{R})}+R^2|F(0)|+R^2|A|\right),
~\forall ~x\in B_{R/2},
\end{equation*}
where $0<\alpha<1$ and $C>0$ are universal constants. Since $\|u\|_{L^{\infty}(B_{R})}\leq K(R^2+1)$, for any fixed $x\in \mathbb{R}^n$, by letting $R\to \infty$, we have
\begin{equation*}
  u(x)=P(x).
\end{equation*}
That is, $u$ is a quadratic polynomial.~\qed~\\

Now, we give the~\\
\noindent\textbf{Proof of \Cref{th1.0}.} Throughout this proof, $C$ always denotes a universal constant. Let $v_i$ ($i\geq 2$) be viscosity solutions of
\begin{equation*}
  \left\{
  \begin{aligned}
    F(D^2v_i)&=0~~\mbox{ in}~~B_i;\\
    v_i&=u~~\mbox{ on}~~\partial B_i.
  \end{aligned}
  \right.
\end{equation*}
For the existence of viscosity solutions, we refer to \cite[Theorem 4.1]{MR1118699} or \cite[Theorem 4.17]{MR4560756}. Then there exist infinite $i$ such that
\begin{equation}\label{e1.3-2}
\max_{x\in \partial B_1}(v_i(x)-u(x))\geq 0
\end{equation}
or
\begin{equation}\label{e1.9-2}
\min_{x\in \partial B_1}(v_i(x)-u(x))\leq 0.
\end{equation}
Let us suppose that \cref{e1.3-2} holds. By the continuity, there exists $x_i\in \partial B_1$ such that
\begin{equation*}
a_i:=\max_{x\in \partial B_1}(v_i(x)-u(x))=v_i(x_i)-u(x_i)\geq 0.
\end{equation*}
Then
\begin{equation}\label{e1.5-2}
  v_i-a_i\leq u~~\mbox{ on }~~\partial B_1~\mbox{and}~\partial B_{i}.
\end{equation}

For $x\in B_i$, let
\begin{equation*}
w_i(x)=v_i(x)-a_i=v_i(x)-v_i(x_i)+u(x_i).
\end{equation*}
Then for any $R\geq 1$ and $i\geq 2R$, with the aid of the interior $C^{0,1}$ estimate for $v_i$, the maximum principle (see \cite[Corollary 3.7]{MR1351007}) and \cref{e1.2-2}, we have
\begin{equation}\label{e1.6-2}
  \begin{aligned}
&\|Dw_i\|_{L^{\infty}(B_{R})}=\|Dv_i\|_{L^{\infty}(B_{R})}\leq \frac{C}{i}\|v_i\|_{L^{\infty}(B_{i})}\leq\frac{C}{i}\|u\|_{L^{\infty}(B_{i})}\leq C, \\
&\|w_i\|_{L^{\infty}(B_{R})}\leq \|v_i-v_i(x_i)\|_{L^{\infty}(B_{R})}+|u(x_i)|
\leq 2R\|Dv_i\|_{L^{\infty}(B_{R})}+1\leq CR.
  \end{aligned}
\end{equation}

By the Arzel\`{a}-Ascoli theorem, there exist $w\in C^{0,1}_{loc}(\mathbb{R}^n)$ and a subsequence of $w_i$ (denoted by $w_i$ again) such that $w_i\rightarrow w$ in $L^{\infty}_{loc}(\mathbb{R}^n)$. From the closedness of viscosity solutions (see \cite[Proposition 2.9]{MR1351007}), $w$ is a viscosity solution of
\begin{equation*}
F(D^2w)=0~~\mbox{ in}~~\mathbb{R}^n.
\end{equation*}
In addition, from \cref{e1.6-2},
\begin{equation}\label{e5.4}
w(x)\leq C(|x|+1),~\forall ~x\in \mathbb{R}^n.
\end{equation}
By \Cref{le3.1}, $w\equiv P$ for some linear polynomial $P$.

From \cref{e1.5-2} and the comparison principle, for any $i\geq 2$,
\begin{equation*}
w_i(x)\leq u(x),~\forall ~x\in B_i\backslash B_1.
\end{equation*}
By letting $i\rightarrow \infty$,
\begin{equation*}
P(x)\leq u(x),~~\forall ~x\in B_1^{c}.
\end{equation*}

If there exist infinite $i$ such that \cref{e1.9-2} holds, by similar arguments as above, we can prove that there exists a linear polynomial $Q$ such that
\begin{equation*}
u\leq Q~~\mbox{ in}~~B_1^c.
\end{equation*}

Next, we assume $u\geq P$ (the case $u\leq Q$ can be proved similarly). Let $v=u-P\geq 0$ and then $v$ is a viscosity solution of \cref{e.F1}. Since $F(0)=0$,
\begin{equation*}
u\in S(\lambda,\Lambda,0)~~\mbox{ in}~\bar B_1^c.
\end{equation*}
In addition, from \cref{e1.2-2} and \cref{e5.4},
\begin{equation*}
\|v\|_{L^{\infty}(B_4\backslash B_1)}\leq \|u\|_{L^{\infty}(B_4\backslash B_1)}+\|P\|_{L^{\infty}(B_4\backslash B_1)}\leq C.
\end{equation*}
Then by \Cref{le2.3}, for any $x\in B_2^c$,
\begin{equation}\label{e5.1}
\left\{
\begin{aligned}
&|v(x)-v_{\infty}|\leq C\|v-v_{\infty}\|_{L^{\infty}(\partial B_1)}E^+(x),~&&\alpha^+>0;\\
&|v(x)|\leq -CE^+(x),~&&\alpha^+\leq 0.
\end{aligned}
\right.
\end{equation}

Finally, from the interior $C^{0,1}$ estimate for $v$, for any $x\in B_4^c$,
\begin{equation}\label{e3.4}
|Dv(x)|\leq
\left\{
\begin{aligned}
&\frac{C}{|x|}\|v-v_{\infty}\|_{L^{\infty}(B(x,|x|/2))}\leq C\|v-v_{\infty}\|_{L^{\infty}(\partial B_1)}|x|^{-1}E^+(x),~&&\alpha^+>0;\\
&\frac{C}{|x|}\|v\|_{L^{\infty}(B(x,|x|/2))}\leq -C|x|^{-1}E^+(x),~&&\alpha^+\leq 0.
\end{aligned}
\right.
\end{equation}
That is, \crefrange{e2.1-2}{e4.1} hold. ~\qed~\\

Now, we give the\\
\noindent\textbf{Proof of \Cref{th1.0-2}.} In the proof of \Cref{th1.0}, we have proved that there exists a linear polynomial $P_0$ such that $u\geq P_0$ or $u\leq P_0$ in $B_1^c$. Hence, $a:=\lim_{x\to \infty} u-P_0$ exists. If $a=+\infty$ or $-\infty$, set $P=P_0$; if $a$ is finite, set $P=P_0+a$. Then by applying \Cref{pr2.3} to $u-P$, we know that $u-P$ must satisfy precisely one of the five alternatives.

If $\alpha^*>0$ and $\tilde{\alpha}^*>0$, then $a$ is finite. By applying \Cref{le2.2-2} (i) (with $v=\Phi,w=u-P$) and (iii) (with $v=-\tilde{\Phi},w=u-P$) respectively, we obtain \cref{e4.2}.~\qed~\\

Next, we prove the asymptotic behavior for solutions with quadratic growth.\\
\noindent\textbf{Proof of \Cref{th1.1}.} Throughout this proof, $C$ always denotes a constant depending only on $n,\lambda,\Lambda, |A|$ and $|F(0)|$. As before, let $v_i$ ($i\geq 2$) be viscosity solutions of
\begin{equation*}
  \left\{
  \begin{aligned}
    F(D^2v_i)&=A~~\mbox{ in}~~B_i;\\
    v_i&=u~~\mbox{ on}~~\partial B_i.
  \end{aligned}
  \right.
\end{equation*}
Without loss of generality, we assume that there exist infinite $i$ such that
\begin{equation}\label{e1.3}
\max_{x\in \partial B_1}(v_i(x)-u(x))\geq 0.
\end{equation}
Take $x_i\in \partial B_1$ such that
\begin{equation*}
a_i:=\max_{x\in \partial B_1}(v_i(x)-u(x))=v_i(x_i)-u(x_i)\geq 0.
\end{equation*}
Then
\begin{equation}\label{e1.5}
  v_i-a_i\leq u~~\mbox{ on }~~\partial B_1~\mbox{and}~\partial B_i.
\end{equation}

For $x\in B_i$, let
\begin{equation*}
w_i(x)=v_i(x)-a_i-Dv_i(x_i)\cdot (x-x_i)=v_i(x)-v_i(x_i)-Dv_i(x_i)\cdot (x-x_i)+u(x_i).
\end{equation*}
Then for any $R\geq 1$ and $i\geq 2R$, with the aid of the interior $C^{1,1}$ estimate for $v_i$, the maximum principle and \cref{e1.2},
\begin{equation}\label{e1.6}
  \begin{aligned}
&\|D^2w_i\|_{L^{\infty}(B_{R})}=\|D^2v_i\|_{L^{\infty}(B_{R})}\leq \frac{C}{i^2}\left(\|v_i\|_{L^{\infty}(B_{i})}+i^2|F(0)|+i^2|A|\right)\\
&\quad\quad\quad\quad\quad\quad
\leq\frac{C}{i^2}\left(\|u\|_{L^{\infty}(B_{i})}+i^2|F(0)|+i^2|A|\right)\leq C, \\
&\|Dw_i\|_{L^{\infty}(B_{R})}=\|Dv_i-Dv_i(x_i)\|_{L^{\infty}(B_{R})}
\leq R\|D^2v_i\|_{L^{\infty}(B_{R})}\leq CR, \\
&\|w_i\|_{L^{\infty}(B_{R})}\leq (R+1)^2\|D^2v_i\|_{L^{\infty}(B_{R})}+|u(x_i)|
\leq CR^2.
  \end{aligned}
\end{equation}

By the Arzel\`{a}-Ascoli theorem, there exists $w\in C(\mathbb{R}^n)$ and a subsequence of $w_i$ (denoted by $w_i$ again) such that $w_i\rightarrow w$ in $L^{\infty}_{loc}(\mathbb{R}^n)$. From the closedness of viscosity solutions, $w$ is a viscosity solution of
\begin{equation*}
F(D^2w)=A~~\mbox{ in}~~\mathbb{R}^n.
\end{equation*}
In addition, from \cref{e1.6},
\begin{equation*}
w(x)\leq C(|x|^2+1),~\forall ~x\in \mathbb{R}^n.
\end{equation*}
By the Liouville theorem in the whole space, $w\equiv P_0$ for some quadratic polynomial $P_0$ and $F(D^2P_0)=A$.

From \cref{e1.5} and the comparison principle, for any $i\geq 2$,
\begin{equation}\label{e1.1}
w_i(x)\leq u(x)-Dv_i(x_i)(x-x_i)~~\mbox{ in}~~B_i\backslash B_1.
\end{equation}
By taking
\begin{equation*}
x=x_i+\tau \frac{Dv_i(x_i)}{|Dv_i(x_i)|}\in \partial B_2,~\tau\in [1,3],
\end{equation*}
in \cref{e1.1}, we have
\begin{equation*}
|Dv_i(x_i)|\leq C,~\forall ~i\geq 2.
\end{equation*}
Then there exist subsequences of $Dv_i(x_i)$ and $x_i$ (denoted by themselves again) such that $Dv_i(x_i)\rightarrow b\in \mathbb{R}^n$ and $x_i\rightarrow x^{\ast}\in \partial B_1$. Thus, by letting $i\rightarrow \infty$ in \cref{e1.1},
\begin{equation*}
P_0(x)\leq u(x)-b\cdot (x-x^{\ast}),~~\forall ~x\in B_1^{c}.
\end{equation*}
Set $P(x)=P_0(x)+b\cdot (x-x^{\ast})$. Then
\begin{equation}\label{e1.4}
u\geq P~~\mbox{ in}~~B_1^c.
\end{equation}

The rest proof is similar to that of \Cref{th1.0}. Let $v=u-P\geq 0$ and then $v$ is a viscosity solution of
\begin{equation*}
G(D^2v):=F(D^2v+D^2P)-A=0~~\mbox{ in}~\bar B_1^c.
\end{equation*}
Note that $G(0)=0$. Hence, $v\in S(\lambda,\Lambda,0)$. By \Cref{le2.3}, \cref{e5.1} holds for $v$. Then by applying the interior $C^{0,1}$ and $C^{1,1}$ estimates for $v$ (since $G$ is also convex), \cref{e3.4} holds and for any $x\in B_4^c$,
\begin{equation*}
|D^2v(x)|\leq
\left\{
\begin{aligned}
&\frac{C}{|x|^2}\|v-v_{\infty}\|_{L^{\infty}(B(x,|x|/2))}\leq C\|v-v_{\infty}\|_{L^{\infty}(\partial B_1)}|x|^{-2}E^+(x),~&&\alpha^+>0;\\
&\frac{C}{|x|^2}\|v\|_{L^{\infty}(B(x,|x|/2))}\leq -C|x|^{-2}E^+(x),~&&\alpha^+\leq 0.
\end{aligned}
\right.
\end{equation*}
That is, \crefrange{e2.1-2-2}{e4.1-2} hold.~\qed
~\\

\noindent\textbf{Data availability statement} Data sharing not applicable to this article as no datasets were generated or analysed during the current study.
%
%
%
~\\

\printbibliography

\end{document}